\documentclass[12pt,reqno]{amsart}
\usepackage[margin=1in]{geometry}
\usepackage{amsmath,amssymb,amsthm,mathtools}
\usepackage{booktabs}
\usepackage[colorlinks=true,linkcolor=blue,citecolor=blue,urlcolor=blue]{hyperref}

\newtheorem{theorem}{Theorem}[section]
\newtheorem{proposition}[theorem]{Proposition}
\newtheorem{lemma}[theorem]{Lemma}

\theoremstyle{definition}
\newtheorem{definition}[theorem]{Definition}
\theoremstyle{remark}
\newtheorem{remark}[theorem]{Remark}
\newtheorem*{example*}{Examples}

\numberwithin{equation}{section}

\newcommand{\mex}{\mathrm{mex}}
\newcommand{\rmex}{r\text{-}\mathrm{mex}}
\newcommand{\Fp}{\mathbb{F}_p}
\newcommand{\W}{\mathcal{W}_{\ell,r}}

\begin{document}

\title{The minimal excludant and the least $r$-gaps of $\ell$-regular partitions}

\author{P. Tejash}
\address{Tejash Priya, Department of Mathematics, National Institute of Technology, Raipur, Chhattisgarh 492010.}
\email{tejash1705@gmail.com, tpriya.phd2026.maths@nitrr.ac.in}

\author{S. Bhowmick}
\address{Sourav Bhowmick, Department of Mathematics, National Institute of Technology, Raipur, Chhattisgarh 492010.}
\email{souravbhowmick578@gmail.com, sbhowmick.phd2025.maths@nitrr.ac.in}

\author{N. K. Meher}
\address{Nabin Kumar Meher, Department of Mathematics, National Institute of Technology, Raipur, Chhattisgarh 492010.}
\email{mehernabin@gmail.com, nkmeher.maths@nitrr.ac.in}

\subjclass[2020]{05A17, 11P81, 11P83}
\keywords{Partitions, $\ell$-regular partitions, minimal excludant, least $r$-gap, generating functions, Jacobi triple product Identity}

\begin{abstract}
The minimal excludant (mex) of an integer partition $\lambda$ is the smallest positive integer that is not a part of $\lambda$. In this paper we study an analogue of this statistic for $\ell$-regular partitions with $\ell\geq2$, that is, partitions in which no part is divisible by $\ell$. Recently, Meher \cite{Meher} defined the minimal excludant for $\ell$-regular partitions. Motivated by Meher \cite{Meher}, for an $\ell$-regular partition we define its least $r$-gap as the smallest positive integer not divisible by $\ell$ that appears fewer than $r$ times as a part, and we obtain the two-variable generating function for the number of $\ell$-regular partitions of $n$ with a given least $r$-gap. Let $\sigma^{\mex}_{\ell}(n)$, and $\sigma^{\mex}_{\ell,r}(n)$ be the sum of the of the mex values over all $\ell$-regular partitions of $n$ and least $r$-gaps over all $\ell$-regular partitions of $n$ respectively. Using Jacobi's triple product identity, we derive explicit expressions for the generating function of $\sigma^{\mex}_{\ell}(n)$, and $\sigma^{\mex}_{\ell,r}(n)$, with separate formulas for even and odd $\ell$. We also show that the sequence $\{\sigma^{\mex}_{\ell,r}(n)\}_{n\ge0}$ is weakly increasing for every $\ell\ge2$ and $r\ge1$. Also, using Ingham's Tauberian theorem, we establish Hardy--Ramanujan type asymptotic formulas for $\sigma^{\mex}_{\ell,r}(n)$ when $\ell=4,5,7$ and $r\ge1$ is arbitrary. Lastly, we prove $\sigma^{\mex}_{\ell,r}$ is odd infinitely often, and we establish a quantitative lower bound, showing that the number of such n up to X grows at least like log log X 
\end{abstract}

\maketitle

\section{Introduction}

A partition of a positive integer $n$ is a non-increasing sequence of positive integers, called parts, whose sum is $n$. Following Andrews and Newman~\cite{AN}, the \emph{minimal excludant} of a partition $\lambda$, denoted $\mex(\lambda)$, is the smallest positive integer that does not appear as a part of $\lambda$. This statistic had earlier been studied by Grabner and Knopfmacher~\cite{GK} under the name ``smallest gap''. Writing
\[
\sigma\mex(n)=\sum_{\lambda\vdash n}\mex(\lambda),
\]
Andrews and Newman~\cite{AN} proved that
\[
\sum_{n\ge0}\sigma\mex(n)\,q^n=(-q;q)_\infty^2,
\]
so that $\sigma\mex(n)$ equals the number of partitions of $n$ into distinct parts using two colors. Since then the minimal excludant and its generalizations have attracted considerable attention; see, for instance, Hopkins, Sellers and Stanton~\cite{HSS}, who connected the minimal excludant to Dyson's crank, and the references therein. For more details about minimal excludant one can see \cite{BS1}, \cite{BS2}, and \cite{KBEM}.

Throughout, we use the standard $q$-Pochhammer notation
\[
(a;q)_\infty=\prod_{k\ge0}(1-aq^k),\qquad f_k:=(q^k;q^k)_\infty,
\]
and we assume $|q|<1$. For an integer $\ell\ge2$, a partition is called \emph{$\ell$-regular} if none of its parts is divisible by $\ell$. Let $\mathcal{P}_\ell(n)$ denote the set of all $\ell$-regular partitions of $n$. The generating function for $\ell$-regular partitions is $f_\ell/f_1$.

Very recently, Meher~\cite{Meher} introduced the minimal excludant of an $\ell$-regular partition. Since multiples of $\ell$ can never occur as parts of an $\ell$-regular partition, he defined, for $\lambda\in\mathcal{P}_\ell(n)$,
\[
\mex_\ell(\lambda):=\min\{m\ge1:\ \ell\nmid m,\ m\notin\lambda\},\qquad
\sigma^{\mex}_{\ell}(n):=\sum_{\lambda\in\mathcal{P}_\ell(n)}\mex_\ell(\lambda).
\]
In \cite{Meher}, the author obtained explicit eta-quotient representations of the generating function of $\sigma^{\mex}_{\ell}(n)$ for $\ell=2$ and $\ell=3$. Also the author established Hardy--Ramanujan type asymptotic formulas for $\sigma^{\mex}_{2}(n)$ and $\sigma^{\mex}_{3}(n)$ and then showed that $\sigma^{\mex}_{2}(n)$ is odd precisely at the generalized pentagonal numbers. More generally, $\sigma^{\mex}_{3}(n)$ is odd infinitely often, with at least $\gg\log\log X$ such $n\le X$. He also suggested extending these results to $\ell\ge4$.

Ballantine and Merca \cite{BM} generalized $\mex(\lambda)$ to $\rmex(\lambda)$, which is the least positive integer that does not occur at least $r$ times as a part of $\lambda$ (so $1\text{-}\mex(\lambda)=\mex(\lambda)$), where $r\in \mathbb{N}.$ For more details about $\rmex(\lambda)$ one can read \cite{BT}.
 Motivated by Ballantine and Merca \cite{BM} and  Meher~\cite{Meher}, for a positive integer $r$ we define the following-

\begin{definition}\label{def:rgap}
For an $\ell$-regular partition $\lambda$, the $\mex_{\ell,r}(\lambda)$ is the smallest positive integer $m$ with $\ell\nmid m$ that appears fewer than $r$ times as a part of $\lambda$. Let $p^{\mex}_{\ell,r}(m,n)$ be the number of $\ell$-regular partitions of $n$ with $r$-gap $m$, and
\[
\sigma^{\mex}_{\ell,r}(n)=\sum_{\lambda}\mex_{\ell,r}(\lambda),
\]
the sum over all $\ell$-regular partitions $\lambda$ of $n$. For $r=1$ these are $p^{\mex}_{\ell}(m,n)$ and $\sigma^{\mex}_{\ell}(n)$.
\end{definition}

Throughout, we write
\[
A_\ell(m):=\sum_{\substack{1\le j<m\\ \ell\nmid j}}j,\qquad
B_\ell(m):=\sum_{\substack{1\le j<m\\ \ell\nmid j}}(r\cdot j)=r\,A_\ell(m).
\]

Our first theorem gives the generating function of $\sigma^{\mex}_{\ell}(n)$ as a finite sum of infinite products. 

\begin{theorem}\label{thm:r1}
Let $\ell\ge2$ be an integer. Then we have two cases

If $\ell$ is even, then
\[
\sum_{n\ge0}\sigma^{\mex}_{\ell}(n)q^n=\frac{f_\ell\,f_{\ell(\ell-1)}}{f_1}\sum_{s=1}^{\ell/2}q^{\frac{s(s-1)}{2}}\Bigl(-q^{\frac{\ell(\ell+2s-3)}{2}};q^{\ell(\ell-1)}\Bigr)_\infty\Bigl(-q^{\frac{\ell(\ell-2s+1)}{2}};q^{\ell(\ell-1)}\Bigr)_\infty.
\]
If $\ell$ is odd, then
\begin{multline*}
\sum_{n\ge0}\sigma^{\mex}_{\ell}(n)q^n=\frac{f_\ell\,f_{\ell(\ell-1)}}{f_1}\Biggl[\sum_{s=1}^{(\ell-1)/2}q^{\frac{s(s-1)}{2}}\Bigl(-q^{\frac{\ell(\ell+2s-3)}{2}};q^{\ell(\ell-1)}\Bigr)_\infty\Bigl(-q^{\frac{\ell(\ell-2s+1)}{2}};q^{\ell(\ell-1)}\Bigr)_\infty\\
+q^{\frac{\ell^2-1}{8}}\frac{f_{2\ell(\ell-1)}^2}{f_{\ell(\ell-1)}^2}\Biggr].
\end{multline*}
\end{theorem}

More generally, for an arbitrary $r\ge1$ we prove the following.

\begin{theorem}\label{thm:general}
Let $\ell\ge2$ and $r\ge1$ be integers. Then we have two cases.

If $\ell$ is even, then
\[
\sum_{n\ge0}\sigma^{\mex}_{\ell,r}(n)q^n=\frac{f_\ell\,f_{r\ell(\ell-1)}}{f_1}\sum_{s=1}^{\ell/2}q^{\frac{rs(s-1)}{2}}\Bigl(-q^{\frac{r\ell(\ell+2s-3)}{2}};q^{r\ell(\ell-1)}\Bigr)_\infty\Bigl(-q^{\frac{r\ell(\ell-2s+1)}{2}};q^{r\ell(\ell-1)}\Bigr)_\infty.
\]
If $\ell$ is odd, then
\begin{multline*}
\sum_{n\ge0}\sigma^{\mex}_{\ell,r}(n)q^n=\frac{f_\ell\,f_{r\ell(\ell-1)}}{f_1}\Biggl[\sum_{s=1}^{(\ell-1)/2}q^{\frac{rs(s-1)}{2}}\Bigl(-q^{\frac{r\ell(\ell+2s-3)}{2}};q^{r\ell(\ell-1)}\Bigr)_\infty\Bigl(-q^{\frac{r\ell(\ell-2s+1)}{2}};q^{r\ell(\ell-1)}\Bigr)_\infty\\
+q^{\frac{r(\ell^2-1)}{8}}\frac{f_{2r\ell(\ell-1)}^2}{f_{r\ell(\ell-1)}^2}\Biggr].
\end{multline*}
\end{theorem}
We establish Hardy--Ramanujan type asymptotic formulas for $\sigma^{\mex}_{\ell,r}(n)$ when $\ell=4,5,7$ and $r\ge1$ is arbitrary, extending the cases $\ell=2,3$ of Meher~\cite[Theorems~1.4 and~1.5]{Meher}. In particular, we have the following results.
\begin{theorem}\label{thm:4}
For every integer $r\ge1$, we have
\[
\sigma^{\mex}_{4,r}(n)\sim\frac{1}{2\sqrt{6rn}}\,e^{\pi\sqrt{n/2}},\qquad n\to\infty.
\]
\end{theorem}

\begin{theorem}\label{thm:5}
For every integer $r\ge1$, we have
\[
\sigma^{\mex}_{5,r}(n)\sim\frac{1}{4\sqrt{2rn}}\,e^{2\pi\sqrt{\frac{2n}{15}}},\qquad n\to\infty.
\]
\end{theorem}

\begin{theorem}\label{thm:7}
For every integer $r\ge1$, we have
\[
\sigma^{\mex}_{7,r}(n)\sim\frac{1}{4\sqrt{3rn}}\,e^{2\pi\sqrt{\frac{n}{7}}},\qquad n\to\infty.
\]
\end{theorem}

Finally, we study the parity of $\sigma^{\mex}_{\ell,r}(n)$. 

\begin{theorem}\label{thm:odd}
For every $\ell\ge2$ and every $r\ge1$ there is a constant $C=C(\ell,r)$ such that for all sufficiently large $X$,
\[
\#\bigl\{0\le n\le X:\ \sigma^{\mex}_{\ell,r}(n)\ \text{is odd}\bigr\}\ \ge\ \frac{\log\log X-C}{\log 2}\ \gg\ \log\log X.
\]
In particular, $\sigma^{\mex}_{\ell,r}(n)$ is odd for infinitely many $n$.
\end{theorem}

The paper is organized as follows. In Section~\ref{sec:prelim} we prove some Propositions and lemmas. In Section~\ref{sec:proofs} we prove Theorem~\ref{thm:r1} to Theorem~\ref{thm:odd}.

\section{Preliminaries}\label{sec:prelim}

\begin{proposition}\label{prop:gf1}
	For $\ell\ge2$, let $f_k=(q^k;q^k)_\infty$ and $A_\ell(m)=\sum_{\substack{1\le j<m\\ \ell\nmid j}}j$. Then
	\[
	\sum_{\substack{m\ge1\\ \ell\nmid m}}\sum_{n\ge0}p^{\mex}_{\ell}(m,n)\,z^mq^n=\frac{f_\ell}{f_1}\sum_{\substack{m\ge1\\ \ell\nmid m}}z^mq^{A_\ell(m)}\bigl(1-q^m\bigr).
	\]
\end{proposition}
\begin{proof}
	See \cite{Meher}.
\end{proof}

\begin{proposition}\label{prop:gf}
For a positive integer $\ell\ge2$ we have
\begin{equation}\label{eq:1.1}
\sum_{\substack{m\ge1\\ \ell\nmid m}}\sum_{n=0}^{\infty}p^{\mex}_{\ell,r}(m,n)\,z^mq^n
=\frac{(q^\ell;q^\ell)_\infty}{(q;q)_\infty}\sum_{\substack{m\ge1\\ \ell\nmid m}}z^mq^{B_\ell(m)}\bigl(1-q^{rm}\bigr).
\end{equation}
\end{proposition}

\begin{proof}[Proof of Proposition~\ref{prop:gf}]
Suppose $\mex_{\ell,r}(\lambda)=m$, $\ell\nmid m$. Then
\begin{enumerate}
\item[(1)] every admissible integer $1\le j<m$, $\ell\nmid j$, must occur at least $r$ times;
\item[(2)] the part $m$ appears at most $r-1$ times;
\item[(3)] every admissible integer greater than $m$ may occur with arbitrary multiplicity.
\end{enumerate}
These three conditions contribute, respectively,
\[
\prod_{\substack{1\le j<m\\ \ell\nmid j}}\frac{q^{rj}}{1-q^j},\qquad
\frac{1-q^{rm}}{1-q^m},\qquad
\prod_{\substack{j>m\\ \ell\nmid j}}\frac{1}{1-q^j}.
\]
The compulsory contribution of the parts below $m$ is
\[
q^{B_\ell(m)}\prod_{\substack{u<m\\ \ell\nmid u}}\frac{1}{1-q^u},
\]
the part $m$ contributes
\[
1+q^m+\cdots+q^{m(r-1)}=\frac{1-q^{rm}}{1-q^m},
\]
and the unrestricted larger admissible parts contribute
\[
\prod_{\substack{u>m\\ \ell\nmid u}}\frac{1}{1-q^u}.
\]
Thus, we obtain
\begin{equation}\label{eq:2.1}
\sum_{n=0}^{\infty}p^{\mex}_{\ell,r}(m,n)q^n=q^{B_\ell(m)}\left(\frac{1-q^{rm}}{1-q^m}\right)\bigl(1-q^m\bigr)\prod_{\substack{u\ge1\\ \ell\nmid u}}\frac{1}{1-q^u},\qquad \ell\nmid m.
\end{equation}
Since
\[
\prod_{\substack{u\ge1\\ \ell\nmid u}}\frac{1}{1-q^u}=\frac{(q^\ell;q^\ell)_\infty}{(q;q)_\infty},
\]
we get
\begin{equation}\label{eq:2.2}
\sum_{n=0}^{\infty}p^{\mex}_{\ell,r}(m,n)q^n=q^{B_\ell(m)}\frac{(q^\ell;q^\ell)_\infty}{(q;q)_\infty}\bigl(1-q^{rm}\bigr),\qquad \ell\nmid m.
\end{equation}
Consequently, multiplying \eqref{eq:2.2} by $z^m$ and summing over $m\ge1$, $\ell\nmid m$, we have
\[
\sum_{\substack{m\ge1\\ \ell\nmid m}}\sum_{n=0}^{\infty}p^{\mex}_{\ell,r}(m,n)\,z^mq^n
=\frac{(q^\ell;q^\ell)_\infty}{(q;q)_\infty}\sum_{\substack{m\ge1\\ \ell\nmid m}}z^mq^{B_\ell(m)}\bigl(1-q^{rm}\bigr),
\]
which is \eqref{eq:1.1}.
\end{proof}

\begin{lemma}[Jacobi triple product Identity]\label{lem:jtp}
For $z\ne0$ and $|Q|<1$,
\[
\sum_{k\in\mathbb{Z}}Q^{k^2}z^k=(-Qz;Q^2)_\infty(-Q/z;Q^2)_\infty(Q^2;Q^2)_\infty.
\]
\end{lemma}
\begin{proof}
	See \cite{Andrews}.
\end{proof}
\begin{lemma}\label{lem:Aformula}
Let $k\ge0$ and $1\le s\le\ell-1$. Then
\begin{equation}\label{eq:2.3}
A_\ell(\ell k+s)=\frac{\ell(\ell-1)}{2}k^2+(s-1)\ell k+\frac{s(s-1)}{2}.
\end{equation}
\end{lemma}

\begin{proof}
The integers $j<\ell k+s$ which are not divisible by $\ell$ are the numbers $\ell i+t$ with $0\le i\le k-1$ and $1\le t\le\ell-1$, together with the numbers $\ell k+t$ with $1\le t\le s-1$. Hence
\begin{align*}
A_\ell(\ell k+s)&=\sum_{i=0}^{k-1}\sum_{t=1}^{\ell-1}(\ell i+t)+\sum_{t=1}^{s-1}(\ell k+t)\\
&=\sum_{i=0}^{k-1}\left(\ell(\ell-1)i+\frac{\ell(\ell-1)}{2}\right)+(s-1)\ell k+\frac{s(s-1)}{2}\\
&=\ell(\ell-1)\frac{k(k-1)}{2}+k\frac{\ell(\ell-1)}{2}+(s-1)\ell k+\frac{s(s-1)}{2}\\
&=\frac{\ell(\ell-1)}{2}k^2+(s-1)\ell k+\frac{s(s-1)}{2}.\qedhere
\end{align*}
\end{proof}

\begin{lemma}\label{lem:reflect}
Let $1\le s\le\ell$. For every integer $k$,
\begin{equation}\label{eq:2.4}
A_\ell\bigl(\ell(-k-1)+s\bigr)=A_\ell\bigl(\ell k+(\ell+1-s)\bigr).
\end{equation}
\end{lemma}

\begin{proof}
Replacing $k$ by $-k-1$ in \eqref{eq:2.3}, we get
\begin{align*}
A_\ell\bigl(\ell(-k-1)+s\bigr)&=\frac{\ell(\ell-1)}{2}(-k-1)^2+(s-1)\ell(-k-1)+\frac{s(s-1)}{2}\\
&=\frac{\ell(\ell-1)}{2}k^2+\Bigl[\ell(\ell-1)-(s-1)\ell\Bigr]k+\Bigl[\frac{\ell(\ell-1)}{2}-(s-1)\ell+\frac{s(s-1)}{2}\Bigr].
\end{align*}
The coefficient of $k$ is $\ell(\ell-s)=\bigl((\ell+1-s)-1\bigr)\ell$. The constant term is
\[
\frac{\ell^2+\ell-2s\ell+s^2-s}{2}=\frac{(\ell-s)(\ell-s+1)}{2}=\frac{(\ell+1-s)\bigl((\ell+1-s)-1\bigr)}{2}.
\]
So the whole expression is \eqref{eq:2.3} with $s$ replaced by $\ell+1-s$.
\end{proof}

\begin{lemma}\label{lem:B}
Let $1\le s\le\ell-1$. Then
\begin{enumerate}
\item[(a)] $B_\ell(\ell k+s)=\dfrac{r\ell(\ell-1)}{2}k^2+r(s-1)\ell k+\dfrac{rs(s-1)}{2}$,
\item[(b)] $B_\ell\bigl(\ell(-k-1)+s\bigr)=B_\ell\bigl(\ell k+(\ell+1-s)\bigr)$ for every integer $k$,
\item[(c)] $B_\ell(\ell k+s)+r(\ell k+s)=B_\ell(\ell k+s+1)$ for $1\le s\le\ell-2$,
\item[(d)] $B_\ell(\ell k+\ell-1)+r(\ell k+\ell-1)=\dfrac{r\ell(\ell-1)}{2}(k+1)^2$.
\end{enumerate}
\end{lemma}

\begin{proof}
Parts (a) and (b) follow from Lemmas~\ref{lem:Aformula} and~\ref{lem:reflect}, since $B_\ell=rA_\ell$. For (c), by (a) we have
\[
B_\ell(\ell k+s)+r(\ell k+s)=\frac{r\ell(\ell-1)}{2}k^2+r\bigl((s+1)-1\bigr)\ell k+\frac{r(s+1)s}{2}=B_\ell(\ell k+s+1).
\]
For (d), by (a) with $s=\ell-1$, we get
\[
B_\ell(\ell k+\ell-1)+r(\ell k+\ell-1)=\frac{r\ell(\ell-1)}{2}k^2+r(\ell-1)\ell k+\frac{r\ell(\ell-1)}{2}=\frac{r\ell(\ell-1)}{2}(k+1)^2.\qedhere
\]
\end{proof}

\begin{lemma}\label{lem:eta}
As $y\to0^+$, we have
\begin{equation}\label{eq:2.5}
(e^{-y};e^{-y})_\infty\sim\sqrt{\frac{2\pi}{y}}\,\exp\Bigl(-\frac{\pi^2}{6y}\Bigr).
\end{equation}
\end{lemma}

\begin{proof}
This follows from the transformation formula for the Dedekind eta function. See~\cite{Apostol}.
\end{proof}

\begin{lemma}[Ingham~\cite{Ingham}]\label{lem:ingham}
Let $B(q):=\sum_{n=0}^{\infty}b(n)q^n$ be a power series with radius of convergence $1$. Assume that $\{b(n)\}$ is a weakly increasing sequence of non-negative real numbers. If there are constants $\beta,\gamma\in\mathbb{R}$ and $\mu>0$ such that
\[
B(e^{-y})\sim\beta\,y^{\gamma}e^{\mu/y},\qquad y\to0^+,
\]
then we have
\[
b(n)\sim\frac{\beta\,\mu^{\frac{2\gamma+1}{4}}}{2\sqrt{\pi}\,n^{\frac{2\gamma+3}{4}}}\,e^{2\sqrt{\mu n}},\qquad n\to\infty.
\]
\end{lemma}

\begin{lemma}\label{thm:monotone}
The sequence $\{\sigma^{\mex}_{\ell,r}(n)\}_{n\ge0}$ is weakly increasing for every $\ell\ge2$ and every $r\ge1$.
\end{lemma}

\begin{proof}[Proof of Theorem~\ref{thm:monotone}]
Fix $\ell\ge2$ and $r\ge1$. Consider again
\[
\varphi:\mathcal{P}_\ell(n)\longrightarrow\mathcal{P}_\ell(n+1),\qquad \varphi(\lambda)=\lambda\cup\{1\}.
\]
As before, $\varphi(\lambda)$ is $\ell$-regular of size $n+1$, and $\varphi$ is injective, since we can recover $\lambda$ from $\varphi(\lambda)$ by removing one copy of $1$.

Let $c$ be the number of times $1$ occurs in $\lambda$; then $1$ occurs $c+1$ times in $\varphi(\lambda)$, and all other parts have the same multiplicity in $\lambda$ and in $\varphi(\lambda)$. We look at three cases.

\noindent\emph{Case 1: $c\ge r$.} Then $1$ occurs at least $r$ times in $\lambda$, and also at least $r+1$ times in $\varphi(\lambda)$. So $1$ is not a candidate for the $r$-mex in either partition, and nothing else changes. Hence
\[
\rmex_\ell(\varphi(\lambda))=\rmex_\ell(\lambda).
\]
\noindent\emph{Case 2: $0\le c\le r-2$} (this case only appears when $r\ge2$). Then $1$ occurs $c$ times in $\lambda$ and $c+1\le r-1$ times in $\varphi(\lambda)$, so it occurs fewer than $r$ times in both. Therefore
\[
\rmex_\ell(\lambda)=1=\rmex_\ell(\varphi(\lambda)).
\]
\noindent\emph{Case 3: $c=r-1$.} Then $1$ occurs fewer than $r$ times in $\lambda$, so $\rmex_\ell(\lambda)=1$. But in $\varphi(\lambda)$ the part $1$ occurs exactly $r$ times, so $1$ is no longer a candidate and
\[
\rmex_\ell(\varphi(\lambda))>1=\rmex_\ell(\lambda).
\]

On combining all three cases, we get for every $\lambda\in\mathcal{P}_\ell(n)$
\[
\rmex_\ell(\varphi(\lambda))\ \ge\ \rmex_\ell(\lambda).
\]
Summing over all $\lambda\in\mathcal{P}_\ell(n)$ and using that $\varphi$ is injective, we get
\[
\sigma^{\mex}_{\ell,r}(n+1)=\sum_{\mu\in\mathcal{P}_\ell(n+1)}\rmex_\ell(\mu)\ \ge\ \sum_{\lambda\in\mathcal{P}_\ell(n)}\rmex_\ell(\varphi(\lambda))\ \ge\ \sum_{\lambda\in\mathcal{P}_\ell(n)}\rmex_\ell(\lambda)=\sigma^{\mex}_{\ell,r}(n).
\]
Hence $\sigma^{\mex}_{\ell,r}(n+1)\ge\sigma^{\mex}_{\ell,r}(n)$ for all $n\ge0$, and the sequence is weakly increasing for every $\ell\ge2$ and $r\ge1$.
\end{proof}

 Now we collect all the necessary results that are needed for the proof of Theorem~\ref{thm:odd}. Let us define
\[
F_\ell(q):=\sum_{n\ge0}\sigma^{\mex}_{\ell,r}(n)q^n.
\]
Also for $s\ge1$ and $k\ge0$ let
\[
E_s(k):=\frac{\ell(\ell-1)}{2}k^2+(s-1)\ell k+\frac{s(s-1)}{2}.
\]
One can see that $E_s(k)=A_\ell(\ell k+s)$ and $B_\ell(\ell k+s)=rE_s(k)$ for $1\le s\le\ell-1$, by Lemmas~\ref{lem:Aformula} and~\ref{lem:B}(a). We also define
\[
M:=\begin{cases}\ell, & \ell\ \text{odd},\\ 2\ell, & \ell\ \text{even},\end{cases}
\qquad\text{and}\qquad
\W:=\{0\}\ \cup\ \bigl\{\,rE_s(k)\bmod\ell\ :\ 2\le s\le\ell-1,\ k\ge0\,\bigr\}\ \subseteq\ \mathbb{Z}/\ell\mathbb{Z}.
\]
\begin{remark}
For every integer $a$, we have
\begin{equation}\label{eq:shiftMX}
	\frac{(a+M)(3(a+M)-1)}{2}-\frac{a(3a-1)}{2}=\frac{M(6a+3M-1)}{2}\equiv0\pmod{\ell}.
\end{equation}
Since for $\ell$ odd $M=\ell$ and $3\ell-1$ is even, and for $\ell$ even $M=2\ell$ gives $\ell(6a+6\ell-1)$. So $\frac{a(3a-1)}{2}\bmod\ell$ \emph{depends only on $a\bmod M$}.
\end{remark} 
Now we introduced a condition, which will play an important role in the proof of Theorem \eqref{thm:odd}.\\
\medskip
\noindent\textbf{Condition (H).} One of the following will holds:
\begin{enumerate}
\item[\textbf{(I)}] there is an integer $c\ge2$ with $\dfrac{c(3c-1)}{2}\bmod\ell\notin\W$.
\item[\textbf{(II)}] $r\ge2$ and there is an integer $c\ge2$ with $\Bigl(\dfrac{c(3c-1)}{2}+1\Bigr)\bmod\ell\notin\W$.
\end{enumerate}

\begin{lemma}\label{lem:Fmod2}
For all $\ell\ge2$ and $r\ge1$, we have
\begin{enumerate}
\item[(a)] $\displaystyle f_1F_\ell(q)=f_\ell\left[\sum_{k\in\mathbb{Z}}q^{\frac{r\ell(\ell-1)}{2}k^2}+\sum_{s=2}^{\ell-1}\sum_{k\ge0}q^{rE_s(k)}\right]$,
\item[(b)] $\displaystyle f_1F_\ell(q)\equiv f_\ell\left[1+\sum_{s=2}^{\ell-1}\sum_{k\ge0}q^{rE_s(k)}\right]\pmod 2$,
\item[(c)] the coefficient of $q^n$ in $f_1F_\ell(q)$ is even whenever $n\bmod\ell\notin\W$.
\end{enumerate}
\end{lemma}

\begin{proof}
(a) Differentiating \eqref{eq:1.1} with respect to $z$ and putting $z=1$, and taking $m=\ell k+s$ with $k\ge0$ and $1\le s\le\ell-1$, we get
\[
F_\ell(q)=\frac{f_\ell}{f_1}\sum_{s=1}^{\ell-1}\sum_{k\ge0}(\ell k+s)\,q^{B_\ell(\ell k+s)}\bigl(1-q^{r(\ell k+s)}\bigr).
\]
By Lemma~\ref{lem:B}(c), for $1\le s\le\ell-2$ the exponent $B_\ell(\ell k+s)+r(\ell k+s)$ equals $B_\ell(\ell k+s+1)$, so for $2\le s\le\ell-1$ the total coefficient of $q^{B_\ell(\ell k+s)}$ is $(\ell k+s)-(\ell k+s-1)=1$. By Lemma~\ref{lem:B}(d), the remaining terms are
\[
\sum_{k\ge0}(\ell k+1)q^{\frac{r\ell(\ell-1)}{2}k^2}-\sum_{k\ge0}(\ell k+\ell-1)q^{\frac{r\ell(\ell-1)}{2}(k+1)^2}=1+2\sum_{k\ge1}q^{\frac{r\ell(\ell-1)}{2}k^2}=\sum_{k\in\mathbb{Z}}q^{\frac{r\ell(\ell-1)}{2}k^2}.
\]
Since $B_\ell(\ell k+s)=rE_s(k)$, this proves (a).\\

(b) In (a), $\sum_{k\in\mathbb{Z}}q^{\frac{r\ell(\ell-1)}{2}k^2}=1+2\sum_{k\ge1}q^{\frac{r\ell(\ell-1)}{2}k^2}\equiv1\pmod2$. Therefore, (b) follows.\\

(c) By Euler's pentagonal number theorem, $f_\ell=\sum_{m\in\mathbb{Z}}(-1)^mq^{\ell m(3m-1)/2}$, so every exponent occurring in $f_\ell$ is a multiple of $\ell$. By (b), every exponent occurring in $f_1F_\ell$ modulo $2$ is therefore congruent modulo $\ell$ either to $0$ or to $rE_s(k)$ for some $2\le s\le\ell-1$ and $k\ge0$. That is, its residue lies in $\W$.
\end{proof}

\begin{lemma}\label{lem:euler}
Let $n$ be a positive integer with $n\bmod\ell\notin\W$. Then
\begin{equation}\label{eq:rec}
\sum_{s=0}^{\infty}\sigma^{\mex}_{\ell,r}\Bigl(n-\frac{s(3s-1)}{2}\Bigr)+\sum_{s=1}^{\infty}\sigma^{\mex}_{\ell,r}\Bigl(n-\frac{s(3s+1)}{2}\Bigr)\equiv0\pmod 2.
\end{equation}
\end{lemma}

\begin{proof}
By Euler's pentagonal number theorem, $f_1\equiv\sum_{s\ge0}q^{s(3s-1)/2}+\sum_{s\ge1}q^{s(3s+1)/2}\pmod 2$. Therefore
\[
\Bigl(\sum_{s\ge0}q^{s(3s-1)/2}+\sum_{s\ge1}q^{s(3s+1)/2}\Bigr)\Bigl(\sum_{n\ge0}\sigma^{\mex}_{\ell,r}(n)q^n\Bigr)\equiv f_1F_\ell(q)\pmod 2.
\]
The coefficient of $q^n$ on the left is the left-hand side of \eqref{eq:rec}, and on the right it is even by Lemma~\ref{lem:Fmod2}(c).
\end{proof}

\begin{lemma}\label{lem:window}
Let $a\ge2$ be an integer and let $c\ge2$ be an integer with $a\equiv c\pmod M$.
\begin{enumerate}
\item[(i)] If (I) holds with this $c$, there is $n\in\bigl[2a-1,\ \frac{a(3a-1)}{2}\bigr]$ with $\sigma^{\mex}_{\ell,r}(n)$ odd.
\item[(ii)] If (II) holds with this $c$, there is $n\in\bigl[2a,\ \frac{a(3a-1)}{2}+1\bigr]$ with $\sigma^{\mex}_{\ell,r}(n)$ odd.
\end{enumerate}
\end{lemma}

\begin{proof}
For $s\ge0$ put
\[
S(s):=\frac{a(3a-1)}{2}-\frac{s(3s-1)}{2},\qquad R(s):=\frac{a(3a-1)}{2}-\frac{s(3s+1)}{2}.
\]
Both $S$ and $R$ are strictly decreasing in $s\ge0$, and
\begin{gather*}
S(0)=\frac{a(3a-1)}{2},\quad S(a-1)=3a-2,\quad S(a)=0,\\
R(1)=\frac{a(3a-1)}{2}-2,\quad R(a-1)=2a-1,\quad R(a)=-a.
\end{gather*}

(i) By \eqref{eq:shiftMX}, $\frac{a(3a-1)}{2}\equiv\frac{c(3c-1)}{2}\pmod\ell$, so $n=\frac{a(3a-1)}{2}$ satisfies the hypothesis of Lemma~\ref{lem:euler}. Suppose, for a contradiction, that $\sigma^{\mex}_{\ell,r}(n)$ is even for every $2a-1\le n\le\frac{a(3a-1)}{2}$. Putting $n=\frac{a(3a-1)}{2}$ in \eqref{eq:rec}, we get
\[
\sum_{s\ge0}\sigma^{\mex}_{\ell,r}(S(s))+\sum_{s\ge1}\sigma^{\mex}_{\ell,r}(R(s))\equiv0\pmod2.
\]
For $s\ge a+1$ we have $S(s)<0$, and for $s\ge a$ we have $R(s)<0$; these terms vanish. Moreover $S(a)=0$ and $\sigma^{\mex}_{\ell,r}(0)=1$. Hence
\begin{equation}\label{eq:win1}
1+\sum_{s=0}^{a-1}\sigma^{\mex}_{\ell,r}(S(s))+\sum_{s=1}^{a-1}\sigma^{\mex}_{\ell,r}(R(s))\equiv0\pmod 2.
\end{equation}
All the arguments $S(s)$ $(0\le s\le a-1)$ and $R(s)$ $(1\le s\le a-1)$ lie in $\bigl[2a-1,\frac{a(3a-1)}{2}\bigr]$, because $S(a-1)=3a-2\ge2a-1$ and $R(a-1)=2a-1$. By assumption every value of $\sigma^{\mex}_{\ell,r}$ in \eqref{eq:win1} is even, so the left-hand side is $\equiv1\pmod2$, a contradiction.

(ii) Here $r\ge2$, so $\sigma^{\mex}_{\ell,r}(1)=1$. By \eqref{eq:shiftMX}, $\frac{a(3a-1)}{2}+1\equiv\frac{c(3c-1)}{2}+1\pmod\ell$, so $n=\frac{a(3a-1)}{2}+1$ satisfies the hypothesis of Lemma~\ref{lem:euler}. Suppose that $\sigma^{\mex}_{\ell,r}(n)$ is even for every $2a\le n\le\frac{a(3a-1)}{2}+1$. At $n=\frac{a(3a-1)}{2}+1$, relation \eqref{eq:rec} reads
\[
\sum_{s\ge0}\sigma^{\mex}_{\ell,r}(S(s)+1)+\sum_{s\ge1}\sigma^{\mex}_{\ell,r}(R(s)+1)\equiv0\pmod 2.
\]
For $s\ge a+1$ we have $S(s)+1\le-3a<0$, and for $s\ge a$ we have $R(s)+1\le1-a<0$. Moreover $S(a)+1=1$ and $\sigma^{\mex}_{\ell,r}(1)=1$. Hence
\begin{equation}\label{eq:win2}
1+\sum_{s=0}^{a-1}\sigma^{\mex}_{\ell,r}(S(s)+1)+\sum_{s=1}^{a-1}\sigma^{\mex}_{\ell,r}(R(s)+1)\equiv0\pmod 2.
\end{equation}
The arguments lie in $\bigl[2a,\frac{a(3a-1)}{2}+1\bigr]$, since $S(a-1)+1=3a-1\ge2a$, $R(a-1)+1=2a$ and $S(0)+1=\frac{a(3a-1)}{2}+1$. Again every value of $\sigma^{\mex}_{\ell,r}$ in \eqref{eq:win2} is even by assumption, and the left-hand side is $\equiv1\pmod2$, a contradiction.
\end{proof}

\begin{lemma}\label{lem:witness}
For every $\ell\ge2$ and $r\ge1$, with Condition (H), we have the followings
\begin{enumerate}
\item[(a)] if $3\mid\ell$, then (I) holds with $c=2$ if $r\not\equiv2\pmod3$, and with $c=4$ if $r\equiv2\pmod3$,
\item[(b)] if $3\nmid\ell$ and a prime $p\ge5$ divides $\ell$, then (I) holds if $3r\not\equiv1\pmod p$, and (II) holds if $3r\equiv1\pmod p$ (in which case $r\ge2$),
\item[(c)] if $\ell=2^f$, then (I) holds with $c=\ell$ if $r$ is odd, and with $c=2\ell+1$ if $r$ is even.
\end{enumerate}
\end{lemma}

\begin{proof}
Every $\ell\ge2$ is divisible by $3$, or has a prime divisor $p\ge5$, or is a power of $2$. Therefore every $\ell$ falls under (a), (b) or (c).

\medskip
\noindent\textbf{(a) $3\mid\ell$.} We have $3\mid\ell$ and $3\mid\frac{\ell(\ell-1)}{2}$ (if $\ell$ is odd because $\ell\mid\frac{\ell(\ell-1)}{2}$; if $\ell$ is even because $3\mid\frac{\ell}{2}$). Hence $E_s(k)\equiv\frac{s(s-1)}{2}\pmod3$. The numbers $\frac{s(s-1)}{2}$ are $\equiv0$ or $1\pmod 3$. Hence every element of $\W$ is in $\{0,\ r\bmod 3\}$. Moreover $\frac{c(3c-1)}{2}\equiv c\pmod3$. If $r\not\equiv2\pmod 3$, we can take $c=2$, since $2\notin\{0,r\}\bmod 3$. If $r\equiv2\pmod3$,  then we take $c=4\equiv1\notin\{0,2\}$. In both cases $\frac{c(3c-1)}{2}\bmod\ell\notin\W$.

\medskip
\noindent\textbf{(b)} A prime $p\ge5$ divides $\ell$. As in (a), $p\mid\frac{\ell(\ell-1)}{2}$ and $p\mid\ell$, so $E_s(k)\equiv\frac{s(s-1)}{2}\pmod p$. Now putting $y=2s-1$, which runs over all of $\Fp$, we get $\frac{s(s-1)}{2}=\frac{y^2-1}{8}$. So modulo $p$ every element of $\W$ lies in $\Bigl\{r\cdot\frac{y^2-1}{8}:\ y\in\Fp\Bigr\}$.
 Now define $Q:=\Bigl\{\frac{x^2-1}{24}:\ x\in\Fp\Bigr\}.$\\
Now substituting $x=(6a-1)^2$, we get $ \frac{x^2-1}{24}= \frac{a(3a-1)}{2}$. Since $6a-1$ runs over $\Fp$ as $a$ does, so $\bigl\{\frac{a(3a-1)}{2}\bmod p\bigr\}=Q$. Note that $|Q|=\frac{p+1}{2}<p$. Now we have three sub-cases.
\begin{itemize}
\item $p\mid r$. Then $3rQ=\{0\}$, and $c=p+1$ gives $\frac{c(3c-1)}{2}\equiv1\notin\{0\}$. Hence Condition (I) holds.
\item $p\nmid r$ and $3r\not\equiv1\pmod p$. We claim that $Q\neq 3rQ$. Suppose $Q=3rQ$ (since $3r$ is a unit, $|3rQ|=|Q|$), and hence the sums of the elements of the two sets are same. Let $\mathcal{S}\subseteq\Fp$ be the set of squares (including $0$). Also $\mathcal{S}$ has $\frac{p+1}{2}$ elements and
\[
\sum_{t\in\mathcal{S}}t=\sum_{y=1}^{(p-1)/2}y^2=\frac{(p-1)p(p+1)}{24}\equiv0\pmod p,
\]
since $p\ge5$. As $Q=\{\frac{t-1}{24}:t\in\mathcal{S}\}$, we get $\sum_{x\in Q}x=\frac{0-\frac{p+1}{2}}{24}\equiv-\frac{1}{48}$ and $\sum_{x\in 3rQ}x\equiv-\frac{3r}{48}$. Thus $3r\equiv1$, a contradiction. Hence some $a$ has $\frac{a(3a-1)}{2}\bmod p\notin 3rQ$, and taking $c\ge2$ with $c\equiv a\pmod p$ gives (I).
\item $3r\equiv1\pmod p$. Then $r\not\equiv1$ (otherwise $3\equiv1$), so $r\ge2$. Here $3rQ=Q$, and $Q+1\not\subseteq Q$, because a non-empty proper subset of $\Fp$ is not invariant under $x\mapsto x+1$. So some $a$ has $\bigl(\frac{a(3a-1)}{2}+1\bigr)\bmod p\notin Q=uQ$, and taking $c\ge2$ with $c\equiv a\pmod p$ gives (II).
\end{itemize}
In every sub-case the relevant residue modulo $p$ lies outside the image of $\W$ modulo $p$. Hence the residue modulo $\ell$ lies outside $\W$.

\medskip
\noindent\textbf{(c) $\ell=2^f$.} Here $\frac{\ell(\ell-1)}{2}k^2=\frac{\ell}{2}(\ell-1)k^2\equiv\frac{\ell}{2}k\pmod\ell$, so $E_s(k)\equiv\frac{s(s-1)}{2}+\frac{\ell}{2}k\pmod\ell$. We claim that the triangular numbers $\frac{j(j+1)}{2}$, $0\le j\le\ell-1$, are pairwise incongruent modulo $\ell$. Indeed, $\frac{i(i+1)}{2}\equiv\frac{j(j+1)}{2}\pmod\ell$ means $2^{f+1}\mid(i-j)(i+j+1)$. The two factors have opposite parity, so one of them is divisible by $2^{f+1}$, but $|i-j|<2^f$ and $0<i+j+1<2^{f+1}$, so $i=j$. Taking $j=s-1$, the numbers $\frac{s(s-1)}{2}$ $(1\le s\le\ell)$ cover every residue modulo $\ell$ exactly once, with $s=1\mapsto0$ and $s=\ell\mapsto\frac{\ell(\ell-1)}{2}\equiv\frac{\ell}{2} \pmod\ell$. Hence $\bigl\{\frac{s(s-1)}{2}:2\le s\le\ell-1\bigr\}=\mathbb{Z}/\ell\mathbb{Z}\setminus\{0,\frac{\ell}{2}\}$.
Now if $k$ is even, then
\[E_s(k)\equiv\frac{s(s-1)}{2}+\frac{\ell}{2}k\equiv \frac{s(s-1)}{2}\pmod\ell, \] and
If $k$ is odd 
\[E_s(k)\equiv\frac{s(s-1)}{2}+\frac{\ell}{2}k\equiv \frac{s(s-1)}{2}+\frac{\ell}{2} \pmod\ell. \]
Therefore $\bigl\{E_s(k):2\le s\le\ell-1\bigr\}=\mathbb{Z}/\ell\mathbb{Z}\setminus\{0,\frac{\ell}{2}\}$.
 Hence,
\[
\W=r\cdot\Bigl(\mathbb{Z}/\ell\mathbb{Z}\setminus\bigl\{\tfrac{\ell}{2}\bigr\}\Bigr).
\]
If $r$ is odd, multiplication by $r$ is a bijection fixing $\frac{\ell}{2}$, so $\W=\mathbb{Z}/\ell\mathbb{Z}\setminus\{\frac{\ell}{2}\}$, and with $c=\ell$ we get $\frac{c(3c-1)}{2}=\frac{\ell}{2}(3\ell-1)\equiv\frac{\ell}{2}\notin\W$.\\ 
If $r$ is even, $\W$ consists of even residues, and with $c=2\ell+1$ we get $\frac{c(3c-1)}{2}=(2\ell+1)(3\ell+1)\equiv1\pmod\ell$, an odd residue.
\end{proof}

\section{Proofs of theorems}\label{sec:proofs}

\begin{proof}[Proof of Theorem~\ref{thm:r1}]
By Proposition~\ref{prop:gf1}, we have
\begin{equation}\label{eq:3.1}
\sum_{\substack{m\ge1\\ \ell\nmid m}}\sum_{n\ge0}p^{\mex}_{\ell}(m,n)\,z^mq^n=\frac{f_\ell}{f_1}\sum_{\substack{m\ge1\\ \ell\nmid m}}z^mq^{A_\ell(m)}\bigl(1-q^m\bigr).
\end{equation}
Note that every $m\in\mathbb{N}$ with $\ell\nmid m$ is written in exactly one way as $m=\ell r+s$ with $r\ge0$ and $1\le s\le\ell-1$. Now, differentiating \eqref{eq:3.1} with respect to $z$ and putting $z=1$, we have
\begin{equation}\label{eq:3.2}
\sum_{n\ge0}\sigma^{\mex}_\ell(n)q^n=\frac{f_\ell}{f_1}\sum_{s=1}^{\ell-1}\sum_{r\ge0}(\ell r+s)\,q^{A_\ell(\ell r+s)}\bigl(1-q^{\ell r+s}\bigr).
\end{equation}
Now
\begin{equation}\label{eq:3.3}
\begin{split}
\sum_{n\ge0}\sigma^{\mex}_\ell(n)q^n=\frac{f_\ell}{f_1}\Biggl[\,&\sum_{r\ge0}(\ell r+1)q^{A_\ell(\ell r+1)}-\sum_{r\ge0}(\ell r+1)q^{A_\ell(\ell r+1)+(\ell r+1)}\\
&+\sum_{r\ge0}(\ell r+2)q^{A_\ell(\ell r+2)}-\sum_{r\ge0}(\ell r+2)q^{A_\ell(\ell r+2)+(\ell r+2)}\\
&+\sum_{r\ge0}(\ell r+3)q^{A_\ell(\ell r+3)}-\sum_{r\ge0}(\ell r+3)q^{A_\ell(\ell r+3)+(\ell r+3)}\\
&\qquad\vdots\\
&+\sum_{r\ge0}(\ell r+\ell-1)q^{A_\ell(\ell r+\ell-1)}-\sum_{r\ge0}(\ell r+\ell-1)q^{A_\ell(\ell r+\ell-1)+(\ell r+\ell-1)}\Biggr].
\end{split}
\end{equation}
For $1\le s\le\ell-2$, using \eqref{eq:2.3}, we have
\begin{align*}
A_\ell(\ell r+s)+(\ell r+s)&=\frac{\ell(\ell-1)}{2}r^2+(s-1)\ell r+\ell r+\frac{s(s-1)}{2}+s\\
&=\frac{\ell(\ell-1)}{2}r^2+\bigl((s+1)-1\bigr)\ell r+\frac{(s+1)s}{2}=A_\ell(\ell r+s+1).
\end{align*}
For $s=\ell-1$, we get
\[
A_\ell(\ell r+\ell-1)+(\ell r+\ell-1)=\frac{\ell(\ell-1)}{2}r^2+(\ell-1)\ell r+\frac{(\ell-1)\ell}{2}=\frac{\ell(\ell-1)}{2}\bigl(r^2+2r+1\bigr)=\frac{\ell(\ell-1)}{2}(r+1)^2.
\]
So in \eqref{eq:3.3} the second sum of the row $s$ is $\sum_{r\ge0}(\ell r+s)q^{A_\ell(\ell r+s+1)}$ for $1\le s\le\ell-2$, and it has the same power of $q$ as the first sum of the row $s+1$. Collecting these terms, \eqref{eq:3.3} becomes
\begin{equation}\label{eq:3.4}
\begin{split}
\sum_{n\ge0}\sigma^{\mex}_\ell(n)q^n&=\frac{f_\ell}{f_1}\Biggl[\sum_{r\ge0}(\ell r+1)q^{A_\ell(\ell r+1)}+\sum_{r\ge0}\bigl((\ell r+2)-(\ell r+1)\bigr)q^{A_\ell(\ell r+2)}\\
&\qquad+\sum_{r\ge0}\bigl((\ell r+3)-(\ell r+2)\bigr)q^{A_\ell(\ell r+3)}+\cdots\\
&\qquad+\sum_{r\ge0}\bigl((\ell r+\ell-1)-(\ell r+\ell-2)\bigr)q^{A_\ell(\ell r+\ell-1)}\\
&\qquad-\sum_{r\ge0}(\ell r+\ell-1)q^{\frac{\ell(\ell-1)}{2}(r+1)^2}\Biggr]\\
&=\frac{f_\ell}{f_1}\Biggl[\sum_{r\ge0}(\ell r+1)q^{\frac{\ell(\ell-1)}{2}r^2}-\sum_{r\ge0}(\ell r+\ell-1)q^{\frac{\ell(\ell-1)}{2}(r+1)^2}\\
&\qquad+\sum_{s=2}^{\ell-1}\sum_{r\ge0}q^{A_\ell(\ell r+s)}\Biggr].
\end{split}
\end{equation}
Replacing $r+1$ by $r$ in the second sum of \eqref{eq:3.4}, we get
\[
\sum_{r\ge0}(\ell r+\ell-1)q^{\frac{\ell(\ell-1)}{2}(r+1)^2}=\sum_{r\ge1}(\ell r-1)q^{\frac{\ell(\ell-1)}{2}r^2},
\]
and so
\begin{align*}
\sum_{r\ge0}(\ell r+1)q^{\frac{\ell(\ell-1)}{2}r^2}-\sum_{r\ge1}(\ell r-1)q^{\frac{\ell(\ell-1)}{2}r^2}&=1+\sum_{r\ge1}\bigl((\ell r+1)-(\ell r-1)\bigr)q^{\frac{\ell(\ell-1)}{2}r^2}\\
&=1+2\sum_{r\ge1}q^{\frac{\ell(\ell-1)}{2}r^2}=\sum_{r\in\mathbb{Z}}q^{\frac{\ell(\ell-1)}{2}r^2}=\sum_{r\in\mathbb{Z}}q^{A_\ell(\ell r+1)}.
\end{align*}
Putting this in \eqref{eq:3.4}, we obtain
\begin{equation}\label{eq:3.5}
\sum_{n\ge0}\sigma^{\mex}_\ell(n)q^n=\frac{f_\ell}{f_1}\left[\sum_{r\in\mathbb{Z}}q^{A_\ell(\ell r+1)}+\sum_{s=2}^{\ell-1}\sum_{r\ge0}q^{A_\ell(\ell r+s)}\right].
\end{equation}

Now we use Lemma~\ref{lem:reflect}. Since $r\ge0$ is the same as $-r-1\le-1$, \eqref{eq:2.4} gives, for $2\le s\le\ell-1$,
\begin{equation}\label{eq:3.6}
\sum_{r\ge0}q^{A_\ell(\ell r+\ell+1-s)}=\sum_{r\ge0}q^{A_\ell(\ell(-r-1)+s)}=\sum_{r\le-1}q^{A_\ell(\ell r+s)}.
\end{equation}
Therefore
\begin{equation}\label{eq:3.7}
\sum_{r\ge0}q^{A_\ell(\ell r+s)}+\sum_{r\ge0}q^{A_\ell(\ell r+\ell+1-s)}=\sum_{r\in\mathbb{Z}}q^{A_\ell(\ell r+s)}.
\end{equation}

\medskip
\noindent\textbf{Case 1: $\ell$ is even.} The second sum in \eqref{eq:3.5} has $\ell-2$ terms, namely $s=2,3,\dots,\ell-1$. Since $\ell+1$ is odd, so the indices $2,3,\dots,\ell-1$ split completely into the pairs $\{2,\ell-1\},\{3,\ell-2\},\dots,\{\frac{\ell}{2},\frac{\ell}{2}+1\}$. Now, writing the second sum of \eqref{eq:3.5} according to these pairs, and using the identity \eqref{eq:3.7} for $s=2,3,\dots,\frac{\ell}{2}$, we get
\begin{align*}
\sum_{s=2}^{\ell-1}\sum_{r\ge0}q^{A_\ell(\ell r+s)}&=\left[\sum_{r\ge0}q^{A_\ell(\ell r+2)}+\sum_{r\ge0}q^{A_\ell(\ell r+\ell-1)}\right]+\left[\sum_{r\ge0}q^{A_\ell(\ell r+3)}+\sum_{r\ge0}q^{A_\ell(\ell r+\ell-2)}\right]\\
&\qquad+\cdots+\left[\sum_{r\ge0}q^{A_\ell(\ell r+\frac{\ell}{2})}+\sum_{r\ge0}q^{A_\ell(\ell r+\frac{\ell}{2}+1)}\right]\\
&=\sum_{r\in\mathbb{Z}}q^{A_\ell(\ell r+2)}+\sum_{r\in\mathbb{Z}}q^{A_\ell(\ell r+3)}+\cdots+\sum_{r\in\mathbb{Z}}q^{A_\ell(\ell r+\frac{\ell}{2})}.
\end{align*}
Therefore, we obtain
\begin{equation}\label{eq:3.8}
\begin{split}
\sum_{n\ge0}\sigma^{\mex}_\ell(n)q^n&=\frac{f_\ell}{f_1}\left[\sum_{r\in\mathbb{Z}}q^{\frac{\ell(\ell-1)}{2}r^2}+\sum_{r\in\mathbb{Z}}q^{A_\ell(\ell r+2)}+\cdots+\sum_{r\in\mathbb{Z}}q^{A_\ell(\ell r+\frac{\ell}{2})}\right]\\
&=\frac{f_\ell}{f_1}\sum_{s=1}^{\ell/2}\sum_{r\in\mathbb{Z}}q^{A_\ell(\ell r+s)}.
\end{split}
\end{equation}
By \eqref{eq:2.3}, we find
\[
\sum_{r\in\mathbb{Z}}q^{A_\ell(\ell r+s)}=q^{\frac{s(s-1)}{2}}\sum_{r\in\mathbb{Z}}\Bigl(q^{\frac{\ell(\ell-1)}{2}}\Bigr)^{r^2}\Bigl(q^{(s-1)\ell}\Bigr)^{r}.
\]
Now, we take $Q=q^{\frac{\ell(\ell-1)}{2}}$ and $z=q^{(s-1)\ell}$. Then $Q^2=q^{\ell(\ell-1)}$, and
\[
Qz=q^{\frac{\ell(\ell-1)}{2}+(s-1)\ell}=q^{\frac{\ell(\ell-1)+2\ell(s-1)}{2}}=q^{\frac{\ell(\ell+2s-3)}{2}},\qquad
\frac{Q}{z}=q^{\frac{\ell(\ell-1)}{2}-(s-1)\ell}=q^{\frac{\ell(\ell-2s+1)}{2}}.
\]
Hence
\begin{equation}\label{eq:3.9}
\sum_{r\in\mathbb{Z}}q^{A_\ell(\ell r+s)}=q^{\frac{s(s-1)}{2}}\Bigl(-q^{\frac{\ell(\ell+2s-3)}{2}};q^{\ell(\ell-1)}\Bigr)_\infty\Bigl(-q^{\frac{\ell(\ell-2s+1)}{2}};q^{\ell(\ell-1)}\Bigr)_\infty f_{\ell(\ell-1)}.
\end{equation}
Putting \eqref{eq:3.9} in \eqref{eq:3.8} gives the formula for even $\ell$.

\medskip
\noindent\textbf{Case 2: $\ell$ is odd.} Again the second sum in \eqref{eq:3.5} has $\ell-2$ terms, namely $s=2,3,\dots,\ell-1$, but now $\ell-2$ is odd. The index $s=\frac{\ell+1}{2}$ satisfies $\ell+1-s=s$, so it is paired with itself, and the remaining indices split into the pairs $\{2,\ell-1\},\{3,\ell-2\},\dots,\{\frac{\ell-1}{2},\frac{\ell+3}{2}\}$. Writing the second sum of \eqref{eq:3.5} according to these pairs, and using the identity \eqref{eq:3.7} for $s=2,3,\dots,\frac{\ell-1}{2}$, we get
\begin{align*}
\sum_{s=2}^{\ell-1}\sum_{r\ge0}q^{A_\ell(\ell r+s)}&=\left[\sum_{r\ge0}q^{A_\ell(\ell r+2)}+\sum_{r\ge0}q^{A_\ell(\ell r+\ell-1)}\right]+\left[\sum_{r\ge0}q^{A_\ell(\ell r+3)}+\sum_{r\ge0}q^{A_\ell(\ell r+\ell-2)}\right]\\
&\qquad+\cdots+\left[\sum_{r\ge0}q^{A_\ell(\ell r+\frac{\ell-1}{2})}+\sum_{r\ge0}q^{A_\ell(\ell r+\frac{\ell+3}{2})}\right]+\sum_{r\ge0}q^{A_\ell(\ell r+\frac{\ell+1}{2})}\\
&=\sum_{r\in\mathbb{Z}}q^{A_\ell(\ell r+2)}+\sum_{r\in\mathbb{Z}}q^{A_\ell(\ell r+3)}+\cdots+\sum_{r\in\mathbb{Z}}q^{A_\ell(\ell r+\frac{\ell-1}{2})}+\sum_{r\ge0}q^{A_\ell(\ell r+\frac{\ell+1}{2})}.
\end{align*}
For the middle term $s=\frac{\ell+1}{2}$ we have $(s-1)\ell=\frac{(\ell-1)\ell}{2}$ and $\frac{s(s-1)}{2}=\frac{(\ell+1)(\ell-1)}{8}=\frac{\ell^2-1}{8}$, so by \eqref{eq:2.3}, we get
\[
A_\ell\Bigl(\ell r+\frac{\ell+1}{2}\Bigr)=\frac{\ell(\ell-1)}{2}r^2+\frac{\ell(\ell-1)}{2}r+\frac{\ell^2-1}{8}=\frac{\ell(\ell-1)}{2}\bigl(r^2+r\bigr)+\frac{\ell^2-1}{8}.
\]
Since $(-r-1)^2+(-r-1)=r^2+r$, the term with $r\ge0$ is equal to the term with $-r-1\le-1$, so the sum over $r\ge0$ is exactly half of the sum over $r\in\mathbb{Z}$, that is,
\[
\sum_{r\ge0}q^{A_\ell(\ell r+\frac{\ell+1}{2})}=\frac12q^{\frac{\ell^2-1}{8}}\sum_{r\in\mathbb{Z}}q^{\frac{\ell(\ell-1)}{2}(r^2+r)}.
\]
Therefore, we obtain
\begin{equation}\label{eq:3.10}
\begin{split}
\sum_{n\ge0}\sigma^{\mex}_\ell(n)q^n&=\frac{f_\ell}{f_1}\Biggl[\sum_{r\in\mathbb{Z}}q^{\frac{\ell(\ell-1)}{2}r^2}+\sum_{r\in\mathbb{Z}}q^{A_\ell(\ell r+2)}+\cdots+\sum_{r\in\mathbb{Z}}q^{A_\ell(\ell r+\frac{\ell-1}{2})}\\
&\qquad\qquad+\frac12q^{\frac{\ell^2-1}{8}}\sum_{r\in\mathbb{Z}}q^{\frac{\ell(\ell-1)}{2}(r^2+r)}\Biggr]\\
&=\frac{f_\ell}{f_1}\left[\sum_{s=1}^{(\ell-1)/2}\sum_{r\in\mathbb{Z}}q^{A_\ell(\ell r+s)}+\frac12q^{\frac{\ell^2-1}{8}}\sum_{r\in\mathbb{Z}}q^{\frac{\ell(\ell-1)}{2}(r^2+r)}\right].
\end{split}
\end{equation}
For $\ell$ odd, \eqref{eq:3.9} takes care of the first sum in \eqref{eq:3.10}, and for the last sum we take $Q=z=q^{\frac{\ell(\ell-1)}{2}}$. Then $Q^{r^2}z^r=q^{\frac{\ell(\ell-1)}{2}(r^2+r)}$, $Qz=q^{\ell(\ell-1)}$, $Q/z=1$, $Q^2=q^{\ell(\ell-1)}$, and
\[
\sum_{r\in\mathbb{Z}}q^{\frac{\ell(\ell-1)}{2}(r^2+r)}=\bigl(-q^{\ell(\ell-1)};q^{\ell(\ell-1)}\bigr)_\infty\bigl(-1;q^{\ell(\ell-1)}\bigr)_\infty f_{\ell(\ell-1)}.
\]
Also $(-1;q^{\ell(\ell-1)})_\infty=2\,(-q^{\ell(\ell-1)};q^{\ell(\ell-1)})_\infty$, and $(-x;x)_\infty=\frac{(x^2;x^2)_\infty}{(x;x)_\infty}$. So with $x=q^{\ell(\ell-1)}$ we get $(-q^{\ell(\ell-1)};q^{\ell(\ell-1)})_\infty=\frac{f_{2\ell(\ell-1)}}{f_{\ell(\ell-1)}}$. Therefore
\begin{equation}\label{eq:3.11}
\frac12\sum_{r\in\mathbb{Z}}q^{\frac{\ell(\ell-1)}{2}(r^2+r)}=\left(\frac{f_{2\ell(\ell-1)}}{f_{\ell(\ell-1)}}\right)^2f_{\ell(\ell-1)}=\frac{f_{2\ell(\ell-1)}^2}{f_{\ell(\ell-1)}}.
\end{equation}
Putting \eqref{eq:3.11} and \eqref{eq:3.9} in \eqref{eq:3.10}, and taking $f_{\ell(\ell-1)}$ out of the bracket, we get the formula for $\ell$ odd.
\end{proof}


From now on $r\ge1$ is arbitrary and, as in Section~\ref{sec:prelim}, summation indices are written $k$.

\begin{proof}[Proof of Theorem~\ref{thm:general}]
By Proposition~\ref{prop:gf}, we have
\begin{equation}\label{eq:4.1}
\sum_{\substack{m\ge1\\ \ell\nmid m}}\sum_{n\ge0}p^{\mex}_{\ell,r}(m,n)\,z^mq^n=\frac{f_\ell}{f_1}\sum_{\substack{m\ge1\\ \ell\nmid m}}z^mq^{B_\ell(m)}\bigl(1-q^{rm}\bigr).
\end{equation}
Every $m\in\mathbb{N}$ with $\ell\nmid m$ is written in exactly one way as $m=\ell k+s$ with $k\ge0$ and $1\le s\le\ell-1$. Differentiating \eqref{eq:4.1} with respect to $z$ and putting $z=1$, we have
\begin{equation}\label{eq:4.2}
\sum_{n\ge0}\sigma^{\mex}_{\ell,r}(n)q^n=\frac{f_\ell}{f_1}\sum_{s=1}^{\ell-1}\sum_{k\ge0}(\ell k+s)\,q^{B_\ell(\ell k+s)}\bigl(1-q^{r(\ell k+s)}\bigr).
\end{equation}
Written out row by row in $s$, we get
\begin{equation}\label{eq:4.3}
\begin{split}
\sum_{n\ge0}\sigma^{\mex}_{\ell,r}(n)q^n=\frac{f_\ell}{f_1}\Biggl[\,&\sum_{k\ge0}(\ell k+1)q^{B_\ell(\ell k+1)}-\sum_{k\ge0}(\ell k+1)q^{B_\ell(\ell k+1)+r(\ell k+1)}\\
&+\sum_{k\ge0}(\ell k+2)q^{B_\ell(\ell k+2)}-\sum_{k\ge0}(\ell k+2)q^{B_\ell(\ell k+2)+r(\ell k+2)}\\
&+\cdots\\
&+\sum_{k\ge0}(\ell k+\ell-1)q^{B_\ell(\ell k+\ell-1)}-\sum_{k\ge0}(\ell k+\ell-1)q^{B_\ell(\ell k+\ell-1)+r(\ell k+\ell-1)}\Biggr].
\end{split}
\end{equation}
By Lemma~\ref{lem:B}(c), for $1\le s\le\ell-2$ the second sum of row $s$ is $\sum_{k\ge0}(\ell k+s)q^{B_\ell(\ell k+s+1)}$, which has the same powers of $q$ as the first sum of row $s+1$. By Lemma~\ref{lem:B}(d), the second sum of the last row is $\sum_{k\ge0}(\ell k+\ell-1)q^{\frac{r\ell(\ell-1)}{2}(k+1)^2}$. Collecting terms, and using $B_\ell(\ell k+1)=\frac{r\ell(\ell-1)}{2}k^2$,
\begin{equation}\label{eq:4.4}
\sum_{n\ge0}\sigma^{\mex}_{\ell,r}(n)q^n=\frac{f_\ell}{f_1}\left[\sum_{k\ge0}(\ell k+1)q^{\frac{r\ell(\ell-1)}{2}k^2}-\sum_{k\ge0}(\ell k+\ell-1)q^{\frac{r\ell(\ell-1)}{2}(k+1)^2}+\sum_{s=2}^{\ell-1}\sum_{k\ge0}q^{B_\ell(\ell k+s)}\right].
\end{equation}
Replacing $k+1$ by $k$ in the second sum of \eqref{eq:4.4}, we get
\[
\sum_{k\ge0}(\ell k+\ell-1)q^{\frac{r\ell(\ell-1)}{2}(k+1)^2}=\sum_{k\ge1}(\ell k-1)q^{\frac{r\ell(\ell-1)}{2}k^2},
\]
and so
\[
\sum_{k\ge0}(\ell k+1)q^{\frac{r\ell(\ell-1)}{2}k^2}-\sum_{k\ge1}(\ell k-1)q^{\frac{r\ell(\ell-1)}{2}k^2}=1+2\sum_{k\ge1}q^{\frac{r\ell(\ell-1)}{2}k^2}=\sum_{k\in\mathbb{Z}}q^{B_\ell(\ell k+1)}.
\]
Putting this in \eqref{eq:4.4}, we obtain
\begin{equation}\label{eq:4.5}
\sum_{n\ge0}\sigma^{\mex}_{\ell,r}(n)q^n=\frac{f_\ell}{f_1}\left[\sum_{k\in\mathbb{Z}}q^{B_\ell(\ell k+1)}+\sum_{s=2}^{\ell-1}\sum_{k\ge0}q^{B_\ell(\ell k+s)}\right].
\end{equation}

Now we use Lemma~\ref{lem:B}(b). Since $k\ge0$ is the same as $-k-1\le-1$, for $2\le s\le\ell-1$ we have
\begin{equation}\label{eq:4.6}
\sum_{k\ge0}q^{B_\ell(\ell k+\ell+1-s)}=\sum_{k\ge0}q^{B_\ell(\ell(-k-1)+s)}=\sum_{k\le-1}q^{B_\ell(\ell k+s)},
\end{equation}
and therefore
\begin{equation}\label{eq:4.7}
\sum_{k\ge0}q^{B_\ell(\ell k+s)}+\sum_{k\ge0}q^{B_\ell(\ell k+\ell+1-s)}=\sum_{k\in\mathbb{Z}}q^{B_\ell(\ell k+s)}.
\end{equation}
By Lemma~\ref{lem:B}(a),
\[
\sum_{k\in\mathbb{Z}}q^{B_\ell(\ell k+s)}=q^{\frac{rs(s-1)}{2}}\sum_{k\in\mathbb{Z}}q^{\frac{r\ell(\ell-1)}{2}k^2}q^{r(s-1)\ell k}.
\]
In Lemma~\ref{lem:jtp} take $Q=q^{\frac{r\ell(\ell-1)}{2}}$ and $z=q^{r(s-1)\ell}$. Then $Q^2=q^{r\ell(\ell-1)}$ and
\[
Qz=q^{\frac{r\ell(\ell-1)}{2}+r(s-1)\ell}=q^{\frac{r\ell(\ell+2s-3)}{2}},\qquad
\frac{Q}{z}=q^{\frac{r\ell(\ell-1)}{2}-r(s-1)\ell}=q^{\frac{r\ell(\ell-2s+1)}{2}}.
\]
Hence
\begin{equation}\label{eq:4.8}
\sum_{k\in\mathbb{Z}}q^{B_\ell(\ell k+s)}=q^{\frac{rs(s-1)}{2}}\Bigl(-q^{\frac{r\ell(\ell+2s-3)}{2}};q^{r\ell(\ell-1)}\Bigr)_\infty\Bigl(-q^{\frac{r\ell(\ell-2s+1)}{2}};q^{r\ell(\ell-1)}\Bigr)_\infty f_{r\ell(\ell-1)}.
\end{equation}

\medskip
\noindent\textbf{Case 1: $\ell$ is even.} The second sum in \eqref{eq:4.5} has $\ell-2$ terms, namely $s=2,3,\dots,\ell-1$. Since $\ell+1$ is odd, these indices split completely into the pairs $\{2,\ell-1\},\{3,\ell-2\},\dots,\{\frac{\ell}{2},\frac{\ell}{2}+1\}$. Writing the second sum of \eqref{eq:4.5} according to these pairs and using \eqref{eq:4.7} for $s=2,3,\dots,\frac{\ell}{2}$, we get
\[
\sum_{s=2}^{\ell-1}\sum_{k\ge0}q^{B_\ell(\ell k+s)}=\sum_{s=2}^{\ell/2}\sum_{k\in\mathbb{Z}}q^{B_\ell(\ell k+s)}.
\]
Therefore
\begin{equation}\label{eq:4.9}
\sum_{n\ge0}\sigma^{\mex}_{\ell,r}(n)q^n=\frac{f_\ell}{f_1}\sum_{s=1}^{\ell/2}\sum_{k\in\mathbb{Z}}q^{B_\ell(\ell k+s)}.
\end{equation}
Putting \eqref{eq:4.8} in \eqref{eq:4.9} gives the formula for even $\ell$.

\medskip
\noindent\textbf{Case 2: $\ell$ is odd.} Again the second sum in \eqref{eq:4.5} has $\ell-2$ terms, namely $s=2,3,\dots,\ell-1$, but now $\ell-2$ is odd. The index $s=\frac{\ell+1}{2}$ satisfies $\ell+1-s=s$, so it is paired with itself, and the remaining indices split into the pairs $\{2,\ell-1\},\{3,\ell-2\},\dots,\{\frac{\ell-1}{2},\frac{\ell+3}{2}\}$. Using \eqref{eq:4.7} for $s=2,3,\dots,\frac{\ell-1}{2}$, we get
\[
\sum_{s=2}^{\ell-1}\sum_{k\ge0}q^{B_\ell(\ell k+s)}=\sum_{s=2}^{(\ell-1)/2}\sum_{k\in\mathbb{Z}}q^{B_\ell(\ell k+s)}+\sum_{k\ge0}q^{B_\ell(\ell k+\frac{\ell+1}{2})}.
\]
For the middle term $s=\frac{\ell+1}{2}$ we have $r(s-1)\ell=\frac{r\ell(\ell-1)}{2}$ and $\frac{rs(s-1)}{2}=\frac{r(\ell+1)(\ell-1)}{8}=\frac{r(\ell^2-1)}{8}$, so by Lemma~\ref{lem:B}(a),
\[
B_\ell\Bigl(\ell k+\frac{\ell+1}{2}\Bigr)=\frac{r\ell(\ell-1)}{2}\bigl(k^2+k\bigr)+\frac{r(\ell^2-1)}{8}.
\]
Since $(-k-1)^2+(-k-1)=k^2+k$, the term with $k\ge0$ equals the term with $-k-1\le-1$, so the sum over $k\ge0$ is exactly half of the sum over $k\in\mathbb{Z}$:
\[
\sum_{k\ge0}q^{B_\ell(\ell k+\frac{\ell+1}{2})}=\frac12q^{\frac{r(\ell^2-1)}{8}}\sum_{k\in\mathbb{Z}}q^{\frac{r\ell(\ell-1)}{2}(k^2+k)}.
\]
Therefore
\begin{equation}\label{eq:4.10}
\sum_{n\ge0}\sigma^{\mex}_{\ell,r}(n)q^n=\frac{f_\ell}{f_1}\left[\sum_{s=1}^{(\ell-1)/2}\sum_{k\in\mathbb{Z}}q^{B_\ell(\ell k+s)}+\frac12q^{\frac{r(\ell^2-1)}{8}}\sum_{k\in\mathbb{Z}}q^{\frac{r\ell(\ell-1)}{2}(k^2+k)}\right].
\end{equation}
The first sum in \eqref{eq:4.10} is given by \eqref{eq:4.8}. For the last sum take $Q=z=q^{\frac{r\ell(\ell-1)}{2}}$ in Lemma~\ref{lem:jtp}. Then $Q^{k^2}z^k=q^{\frac{r\ell(\ell-1)}{2}(k^2+k)}$, $Qz=Q^2=q^{r\ell(\ell-1)}$ and $Q/z=1$, so
\[
\sum_{k\in\mathbb{Z}}q^{\frac{r\ell(\ell-1)}{2}(k^2+k)}=\bigl(-q^{r\ell(\ell-1)};q^{r\ell(\ell-1)}\bigr)_\infty\bigl(-1;q^{r\ell(\ell-1)}\bigr)_\infty f_{r\ell(\ell-1)}.
\]
Also $(-1;q^{r\ell(\ell-1)})_\infty=2(-q^{r\ell(\ell-1)};q^{r\ell(\ell-1)})_\infty$, and $(-x;x)_\infty=\frac{(x^2;x^2)_\infty}{(x;x)_\infty}$. So with $x=q^{r\ell(\ell-1)}$ we get $(-q^{r\ell(\ell-1)};q^{r\ell(\ell-1)})_\infty=\frac{f_{2r\ell(\ell-1)}}{f_{r\ell(\ell-1)}}$. Therefore
\begin{equation}\label{eq:4.11}
\frac12\sum_{k\in\mathbb{Z}}q^{\frac{r\ell(\ell-1)}{2}(k^2+k)}=\left(\frac{f_{2r\ell(\ell-1)}}{f_{r\ell(\ell-1)}}\right)^2f_{r\ell(\ell-1)}=\frac{f_{2r\ell(\ell-1)}^2}{f_{r\ell(\ell-1)}}.
\end{equation}
Putting \eqref{eq:4.11} and \eqref{eq:4.8} in \eqref{eq:4.10}, and taking $f_{r\ell(\ell-1)}$ out of the bracket, we get the formula for odd $\ell$.
\end{proof}

\begin{remark}\label{rem:r1}
For $r=1$ we have $B_\ell=A_\ell$, and Theorem~\ref{thm:general} reduces to Theorem~\ref{thm:r1}, the generating function for $\sigma^{\mex}_\ell(n)$.
\end{remark}

\begin{proof}[Proof of Theorem~\ref{thm:4}]
We have
\begin{align*}
	\sum_{n\ge0}\sigma^{\mex}_{4,r}(n)q^n
	&=\frac{(q^4;q^4)_\infty(q^{12r};q^{12r})_\infty^5}{(q;q)_\infty(q^{6r};q^{6r})_\infty^2(q^{24r};q^{24r})_\infty^2}\\
	&\quad+q^r\frac{(q^4;q^4)_\infty(q^{4r};q^{4r})_\infty^2(q^{6r};q^{6r})_\infty(q^{24r};q^{24r})_\infty}{(q;q)_\infty(q^{2r};q^{2r})_\infty(q^{8r};q^{8r})_\infty(q^{12r};q^{12r})_\infty}.
\end{align*}
Let
\[
T_1(q):=\frac{(q^4;q^4)_\infty(q^{12r};q^{12r})_\infty^5}{(q;q)_\infty(q^{6r};q^{6r})_\infty^2(q^{24r};q^{24r})_\infty^2}.
\]
Using \eqref{eq:2.5}, we get
\[
(e^{-4y};e^{-4y})_\infty\sim\sqrt{\frac{\pi}{2y}}\,e^{-\frac{\pi^2}{24y}},\qquad
(e^{-12ry};e^{-12ry})_\infty^5\sim\Bigl(\frac{\pi}{6ry}\Bigr)^{5/2}e^{-\frac{5\pi^2}{72ry}},
\]
and
\[
(e^{-y};e^{-y})_\infty(e^{-6ry};e^{-6ry})_\infty^2(e^{-24ry};e^{-24ry})_\infty^2\sim\sqrt{\frac{2\pi}{y}}\,\frac{\pi}{3ry}\,\frac{\pi}{12ry}\,e^{-\frac{\pi^2}{6y}-\frac{\pi^2}{18ry}-\frac{\pi^2}{72ry}}.
\]
Since $-\frac{5}{72}+\frac{1}{18}+\frac{1}{72}=0$, therefore,
\[
T_1(e^{-y})\sim\frac12\sqrt{\frac{\pi}{6r}}\;y^{-1/2}e^{\frac{\pi^2}{8y}}.
\]
Similarly, let
\[
T_2(q):=q^r\frac{(q^4;q^4)_\infty(q^{4r};q^{4r})_\infty^2(q^{6r};q^{6r})_\infty(q^{24r};q^{24r})_\infty}{(q;q)_\infty(q^{2r};q^{2r})_\infty(q^{8r};q^{8r})_\infty(q^{12r};q^{12r})_\infty}.
\]
Using \eqref{eq:2.5}, we get
\[
(e^{-4ry};e^{-4ry})_\infty^2(e^{-6ry};e^{-6ry})_\infty(e^{-24ry};e^{-24ry})_\infty\sim\frac{\pi}{2ry}\sqrt{\frac{\pi}{3ry}}\sqrt{\frac{\pi}{12ry}}\;e^{-\frac{\pi^2}{12ry}-\frac{\pi^2}{36ry}-\frac{\pi^2}{144ry}},
\]
and
\[
(e^{-2ry};e^{-2ry})_\infty(e^{-8ry};e^{-8ry})_\infty(e^{-12ry};e^{-12ry})_\infty\sim\sqrt{\frac{\pi}{ry}}\sqrt{\frac{\pi}{4ry}}\sqrt{\frac{\pi}{6ry}}\;e^{-\frac{\pi^2}{12ry}-\frac{\pi^2}{48ry}-\frac{\pi^2}{72ry}}.
\]
Both exponents are equal to $-\frac{17\pi^2}{144ry}$. As $y\to0^+$, we get $q=e^{-y}\to1$ and hence $q^r\sim1$. Thus
\[
T_2(e^{-y})\sim\frac12\sqrt{\frac{\pi}{6r}}\;y^{-1/2}e^{\frac{\pi^2}{8y}}.
\]
Note that both terms have the same exponential growth. Consequently,
\[
T_1(e^{-y})+T_2(e^{-y})=\sum_{n\ge0}\sigma^{\mex}_{4,r}(n)e^{-ny}\sim\Bigl(\frac12\sqrt{\frac{\pi}{6r}}+\frac12\sqrt{\frac{\pi}{6r}}\Bigr)y^{-1/2}e^{\frac{\pi^2}{8y}}=\sqrt{\frac{\pi}{6r}}\;y^{-1/2}e^{\frac{\pi^2}{8y}}.
\]
From Theorem~\ref{thm:monotone}, we know that the sequence $\{\sigma^{\mex}_{4,r}(n)\}$ is weakly increasing. Thus, by applying Lemma~\ref{lem:ingham} with parameters $\beta=\sqrt{\frac{\pi}{6r}}$, $\gamma=-\frac12$, $\mu=\frac{\pi^2}{8}$, we obtain
\[
\sigma^{\mex}_{4,r}(n)\sim\frac{1}{2\sqrt{6rn}}\,e^{\pi\sqrt{n/2}},\qquad n\to\infty.\qedhere
\]
\end{proof}

\begin{proof}[Proof of Theorem~\ref{thm:5}]
We have
\begin{multline*}
\sum_{n\ge0}\sigma^{\mex}_{5,r}(n)q^n=\frac{(q^5;q^5)_\infty(q^{20r};q^{20r})_\infty^5}{(q;q)_\infty(q^{10r};q^{10r})_\infty^2(q^{40r};q^{40r})_\infty^2}\\+q^r\frac{(q^5;q^5)_\infty(q^{10r};q^{10r})_\infty^2}{(q;q)_\infty(q^{5r};q^{5r})_\infty}+q^{3r}\frac{(q^5;q^5)_\infty(q^{40r};q^{40r})_\infty^2}{(q;q)_\infty(q^{20r};q^{20r})_\infty}.
\end{multline*}
Using \eqref{eq:2.5}, we first note that
\[
\frac{(e^{-5y};e^{-5y})_\infty}{(e^{-y};e^{-y})_\infty}\sim\frac{1}{\sqrt5}\,e^{\frac{2\pi^2}{15y}}.
\]
Let
\[
T_1(q):=\frac{(q^5;q^5)_\infty(q^{20r};q^{20r})_\infty^5}{(q;q)_\infty(q^{10r};q^{10r})_\infty^2(q^{40r};q^{40r})_\infty^2}.
\]
Using \eqref{eq:2.5}, we get
\[
(e^{-20ry};e^{-20ry})_\infty^5\sim\Bigl(\frac{\pi}{10ry}\Bigr)^{5/2}e^{-\frac{\pi^2}{24ry}},
\]
and
\[
(e^{-10ry};e^{-10ry})_\infty^2(e^{-40ry};e^{-40ry})_\infty^2\sim\frac{\pi}{5ry}\,\frac{\pi}{20ry}\,e^{-\frac{\pi^2}{30ry}-\frac{\pi^2}{120ry}}.
\]
Since $-\frac{1}{24}+\frac{1}{30}+\frac{1}{120}=0$, therefore,
\[
T_1(e^{-y})\sim\frac15\sqrt{\frac{\pi}{2r}}\;y^{-1/2}e^{\frac{2\pi^2}{15y}}.
\]
Similarly, let
\[
T_2(q):=q^r\frac{(q^5;q^5)_\infty(q^{10r};q^{10r})_\infty^2}{(q;q)_\infty(q^{5r};q^{5r})_\infty}.
\]
Using \eqref{eq:2.5}, we get
\[
(e^{-10ry};e^{-10ry})_\infty^2\sim\frac{\pi}{5ry}\,e^{-\frac{\pi^2}{30ry}},\qquad
(e^{-5ry};e^{-5ry})_\infty\sim\sqrt{\frac{2\pi}{5ry}}\,e^{-\frac{\pi^2}{30ry}}.
\]
As $y\to0^+$, we get $q^r\sim1$. Thus
\[
T_2(e^{-y})\sim\frac15\sqrt{\frac{\pi}{2r}}\;y^{-1/2}e^{\frac{2\pi^2}{15y}}.
\]
Finally, let
\[
T_3(q):=q^{3r}\frac{(q^5;q^5)_\infty(q^{40r};q^{40r})_\infty^2}{(q;q)_\infty(q^{20r};q^{20r})_\infty}.
\]
Using \eqref{eq:2.5}, we get
\[
(e^{-40ry};e^{-40ry})_\infty^2\sim\frac{\pi}{20ry}\,e^{-\frac{\pi^2}{120ry}},\qquad
(e^{-20ry};e^{-20ry})_\infty\sim\sqrt{\frac{\pi}{10ry}}\,e^{-\frac{\pi^2}{120ry}}.
\]
As $y\to0^+$, we get $q^{3r}\sim1$. Thus
\[
T_3(e^{-y})\sim\frac{1}{10}\sqrt{\frac{\pi}{2r}}\;y^{-1/2}e^{\frac{2\pi^2}{15y}}.
\]
Note that all three terms have the same exponential growth. Consequently,
\[
\sum_{n\ge0}\sigma^{\mex}_{5,r}(n)e^{-ny}\sim\Bigl(\frac15+\frac15+\frac1{10}\Bigr)\sqrt{\frac{\pi}{2r}}\;y^{-1/2}e^{\frac{2\pi^2}{15y}}=\frac12\sqrt{\frac{\pi}{2r}}\;y^{-1/2}e^{\frac{2\pi^2}{15y}}.
\]
From Theorem~\ref{thm:monotone}, we know that the sequence $\{\sigma^{\mex}_{5,r}(n)\}$ is weakly increasing. Thus, by applying Lemma~\ref{lem:ingham} with parameters $\beta=\frac12\sqrt{\frac{\pi}{2r}}$, $\gamma=-\frac12$, $\mu=\frac{2\pi^2}{15}$, we obtain
\[
\sigma^{\mex}_{5,r}(n)\sim\frac{1}{4\sqrt{2rn}}\,e^{2\pi\sqrt{\frac{2n}{15}}},\qquad n\to\infty.\qedhere
\]
\end{proof}

\begin{proof}[Proof of Theorem~\ref{thm:7}]
We have
\begin{align*}
\sum_{n\ge0}\sigma^{\mex}_{7,r}(n)q^n&=\frac{(q^7;q^7)_\infty(q^{42r};q^{42r})_\infty^5}{(q;q)_\infty(q^{21r};q^{21r})_\infty^2(q^{84r};q^{84r})_\infty^2}+q^r\frac{(q^7;q^7)_\infty(q^{28r};q^{28r})_\infty(q^{42r};q^{42r})_\infty^2}{(q;q)_\infty(q^{14r};q^{14r})_\infty(q^{84r};q^{84r})_\infty}\\
&\quad+q^{3r}\frac{(q^7;q^7)_\infty(q^{14r};q^{14r})_\infty^2(q^{21r};q^{21r})_\infty(q^{84r};q^{84r})_\infty}{(q;q)_\infty(q^{7r};q^{7r})_\infty(q^{28r};q^{28r})_\infty(q^{42r};q^{42r})_\infty}+q^{6r}\frac{(q^7;q^7)_\infty(q^{84r};q^{84r})_\infty^2}{(q;q)_\infty(q^{42r};q^{42r})_\infty}.
\end{align*}
Using \eqref{eq:2.5}, we first note that
\[
\frac{(e^{-7y};e^{-7y})_\infty}{(e^{-y};e^{-y})_\infty}\sim\frac{1}{\sqrt7}\,e^{\frac{\pi^2}{7y}}.
\]
Let
\[
T_1(q):=\frac{(q^7;q^7)_\infty(q^{42r};q^{42r})_\infty^5}{(q;q)_\infty(q^{21r};q^{21r})_\infty^2(q^{84r};q^{84r})_\infty^2}.
\]
Using \eqref{eq:2.5}, we get
\[
(e^{-42ry};e^{-42ry})_\infty^5\sim\Bigl(\frac{\pi}{21ry}\Bigr)^{5/2}e^{-\frac{5\pi^2}{252ry}},
\]
and
\[
(e^{-21ry};e^{-21ry})_\infty^2(e^{-84ry};e^{-84ry})_\infty^2\sim\frac{2\pi}{21ry}\,\frac{\pi}{42ry}\,e^{-\frac{\pi^2}{63ry}-\frac{\pi^2}{252ry}}.
\]
Since $-\frac{5}{252}+\frac{1}{63}+\frac{1}{252}=0$, therefore,
\[
T_1(e^{-y})\sim\frac17\sqrt{\frac{\pi}{3r}}\;y^{-1/2}e^{\frac{\pi^2}{7y}}.
\]
Similarly, let
\[
T_2(q):=q^r\frac{(q^7;q^7)_\infty(q^{28r};q^{28r})_\infty(q^{42r};q^{42r})_\infty^2}{(q;q)_\infty(q^{14r};q^{14r})_\infty(q^{84r};q^{84r})_\infty}.
\]
Using \eqref{eq:2.5}, we get
\[
(e^{-28ry};e^{-28ry})_\infty(e^{-42ry};e^{-42ry})_\infty^2\sim\sqrt{\frac{\pi}{14ry}}\,\frac{\pi}{21ry}\,e^{-\frac{\pi^2}{168ry}-\frac{\pi^2}{126ry}},
\]
and
\[
(e^{-14ry};e^{-14ry})_\infty(e^{-84ry};e^{-84ry})_\infty\sim\sqrt{\frac{\pi}{7ry}}\sqrt{\frac{\pi}{42ry}}\;e^{-\frac{\pi^2}{84ry}-\frac{\pi^2}{504ry}}.
\]
Both exponents are equal to $-\frac{7\pi^2}{504ry}$. As $y\to0^+$, we get $q^r\sim1$. Thus
\[
T_2(e^{-y})\sim\frac17\sqrt{\frac{\pi}{3r}}\;y^{-1/2}e^{\frac{\pi^2}{7y}}.
\]
Next, let
\[
T_3(q):=q^{3r}\frac{(q^7;q^7)_\infty(q^{14r};q^{14r})_\infty^2(q^{21r};q^{21r})_\infty(q^{84r};q^{84r})_\infty}{(q;q)_\infty(q^{7r};q^{7r})_\infty(q^{28r};q^{28r})_\infty(q^{42r};q^{42r})_\infty}.
\]
Using \eqref{eq:2.5}, we get
\[
(e^{-14ry};e^{-14ry})_\infty^2(e^{-21ry};e^{-21ry})_\infty(e^{-84ry};e^{-84ry})_\infty\sim\frac{\pi}{7ry}\sqrt{\frac{2\pi}{21ry}}\sqrt{\frac{\pi}{42ry}}\;e^{-\frac{\pi^2}{42ry}-\frac{\pi^2}{126ry}-\frac{\pi^2}{504ry}},
\]
and
\[
(e^{-7ry};e^{-7ry})_\infty(e^{-28ry};e^{-28ry})_\infty(e^{-42ry};e^{-42ry})_\infty\sim\sqrt{\frac{2\pi}{7ry}}\sqrt{\frac{\pi}{14ry}}\sqrt{\frac{\pi}{21ry}}\;e^{-\frac{\pi^2}{42ry}-\frac{\pi^2}{168ry}-\frac{\pi^2}{252ry}}.
\]
Both exponents are equal to $-\frac{17\pi^2}{504ry}$. As $y\to0^+$, we get $q^{3r}\sim1$. Thus
\[
T_3(e^{-y})\sim\frac17\sqrt{\frac{\pi}{3r}}\;y^{-1/2}e^{\frac{\pi^2}{7y}}.
\]
Finally, let
\[
T_4(q):=q^{6r}\frac{(q^7;q^7)_\infty(q^{84r};q^{84r})_\infty^2}{(q;q)_\infty(q^{42r};q^{42r})_\infty}.
\]
Using \eqref{eq:2.5}, we get
\[
(e^{-84ry};e^{-84ry})_\infty^2\sim\frac{\pi}{42ry}\,e^{-\frac{\pi^2}{252ry}},\qquad
(e^{-42ry};e^{-42ry})_\infty\sim\sqrt{\frac{\pi}{21ry}}\,e^{-\frac{\pi^2}{252ry}}.
\]
As $y\to0^+$, we get $q^{6r}\sim1$. Thus
\[
T_4(e^{-y})\sim\frac{1}{14}\sqrt{\frac{\pi}{3r}}\;y^{-1/2}e^{\frac{\pi^2}{7y}}.
\]
Note that all four terms have the same exponential growth. Consequently,
\[
\sum_{n\ge0}\sigma^{\mex}_{7,r}(n)e^{-ny}\sim\Bigl(\frac17+\frac17+\frac17+\frac1{14}\Bigr)\sqrt{\frac{\pi}{3r}}\;y^{-1/2}e^{\frac{\pi^2}{7y}}=\frac12\sqrt{\frac{\pi}{3r}}\;y^{-1/2}e^{\frac{\pi^2}{7y}}.
\]
From Theorem~\ref{thm:monotone}, we know that the sequence $\{\sigma^{\mex}_{7,r}(n)\}$ is weakly increasing. Thus, by applying Lemma~\ref{lem:ingham} with parameters $\beta=\frac12\sqrt{\frac{\pi}{3r}}$, $\gamma=-\frac12$, $\mu=\frac{\pi^2}{7}$, we obtain
\[
\sigma^{\mex}_{7,r}(n)\sim\frac{1}{4\sqrt{3rn}}\,e^{2\pi\sqrt{\frac{n}{7}}},\qquad n\to\infty.\qedhere
\]
\end{proof}

\begin{proof}[Proof of Theorem~\ref{thm:odd}]
By Lemma~\ref{lem:witness}, Condition (H) holds; fix $c$ and the case (I) or (II) accordingly. Put $a_0=c$ and
\[
a_{k+1}=M\cdot\frac{a_k(3a_k-1)}{2}+c\qquad(k\ge0).
\]
Since $\frac{a_k(3a_k-1)}{2}$ is an integer, $a_k\equiv c\pmod M$ for all $k$, and $a_k\ge c\ge2$. By Lemma~\ref{lem:window} with $a=a_k$ there is $n_k$ with $\sigma^{\mex}_{\ell,r}(n_k)$ odd and
\begin{equation}\label{eq:5.1}
2a_k-1\le n_k\le\frac{a_k(3a_k-1)}{2}\quad\text{,}\qquad 2a_k\le n_k\le\frac{a_k(3a_k-1)}{2}+1\quad\text{.}
\end{equation}
In both cases $2a_k-1\le n_k\le\frac{a_k(3a_k-1)}{2}+1$.

\emph{The $n_k$ are distinct.} We have $2a_{k+1}-1=M\,a_k(3a_k-1)+2c-1>\frac{a_k(3a_k-1)}{2}+1$, so the windows in \eqref{eq:5.1} are pairwise disjoint and $n_0<n_1<\cdots$.

Let $\nu$ be the largest integer with $\frac{a_\nu(3a_\nu-1)}{2}+1\le X$. Then $n_0,\dots,n_\nu\le X$, so the number of $n\le X$ with $\sigma^{\mex}_{\ell,r}(n)$ odd is at least $\nu+1$.

Since $c\le a_k\le\frac{M}{2}a_k$ (as $M\ge3$),
\[
a_{k+1}=\frac{3M}{2}a_k^2-\frac{M}{2}a_k+c\le\frac{3M}{2}a_k^2<2Ma_k^2,\qquad\text{i.e.}\qquad 2Ma_{k+1}<(2Ma_k)^2.
\]
Iterating, $2Ma_k\le(2Mc)^{2^k}$, so $a_k<(2Mc)^{2^k}$ and $\log\log a_k<k\log2+\log\log(2Mc)$. By the maximality of $\nu$, $X<\frac{a_{\nu+1}(3a_{\nu+1}-1)}{2}+1<a_{\nu+2}$. Therefore
\[
\log\log X<\log\log a_{\nu+2}<\nu\log2+C,\qquad C=2\log2+\log\log(2Mc),
\]
so $\nu>\frac{\log\log X-C}{\log2}$, which proves the theorem.
\end{proof}


\end{document}